\documentclass[a4paper]{amsart}
\usepackage{graphicx}
\usepackage{amsmath}
\usepackage{amssymb}
\usepackage{oldgerm}
\usepackage[all]{xy}

\newcommand{\Hom}{\operatorname{Hom}\nolimits}
\renewcommand{\Im}{\operatorname{Im}\nolimits}

\newcommand{\Ker}{\operatorname{Ker}\nolimits}

\newcommand{\Pd}{\operatorname{pd}\nolimits}
\newcommand{\gldim}{\operatorname{gl.dim}\nolimits}

\newcommand{\Tor}{\operatorname{Tor}\nolimits}
\newcommand{\Ext}{\operatorname{Ext}\nolimits}

\newcommand{\HH}{\operatorname{HH}\nolimits}

\newcommand{\ra}{\operatorname{\mathfrak{r}}\nolimits}
\newcommand{\La}{\operatorname{\Lambda}\nolimits}

\newcommand{\op}{\operatorname{op}\nolimits}

\newcommand{\e}{\operatorname{e}\nolimits}
\newcommand{\Lae}{\operatorname{\Lambda^{\e}}\nolimits}

\theoremstyle{definition}

\theoremstyle{definition}

\theoremstyle{definition}

\theoremstyle{definition}

\theoremstyle{definition}

\theoremstyle{definition}

\theoremstyle{remark}
\newtheorem*{remark}{Remark}
\theoremstyle{definition}

\theoremstyle{definition}

\theoremstyle{definition}

\begin{document}
\title{On the Hochschild (co)homology of quantum exterior algebras}
\author{Petter Andreas Bergh}
\address{Petter Andreas Bergh \\ Institutt for matematiske fag \\
NTNU \\ N-7491 Trondheim \\ Norway} \email{bergh@math.ntnu.no}
%\date{\today}

\subjclass[2000]{16E05, 16E40, 16P10}

\maketitle

\begin{abstract}
We compute the Hochschild cohomology and homology of the algebra
$\La = k \langle x,y \rangle / ( x^2, xy+qyx, y^2 )$ with
coefficients in ${_1 \Lambda_{\psi}}$ for every degree preserving
$k$-algebra automorphism $\psi \colon \La \to \La$. As a result we
obtain several interesting examples of the homological behavior of
$\La$ as a bimodule.
\end{abstract}

\section*{Introduction}

Throughout this paper, let $k$ be a field and $q \in k$ a nonzero
element which is \emph{not} a root of unity. Denote by $\La$ the
$k$-algebra
$$\La = k \langle x,y \rangle / ( x^2, xy+qyx, y^2 ),$$
and by $\Lae$ its enveloping algebra $\La^{\op} \otimes_k \La$.
All modules considered are assumed to be right modules.

During the last years, this $4$-dimensional graded Koszul algebra,
whose module category was classified in \cite{Schulz}, has provided
several examples (or rather counterexamples) giving negative answers
to homological conjectures and questions. Among these are the
conjecture of Auslander on local $\Ext$-limitations (see \cite[page
815]{Auslander}, \cite{Jorgensen} and \cite{Smalo}) and the question
of Happel on the relation between the global dimension and the
vanishing of the Hochschild cohomology (see \cite{Happel} and
\cite{Buchweitz}).

We shall study the Hochschild cohomology and homology of $\La$. More
precisely, for every degree preserving $k$-algebra automorphism
$\psi \colon \La \to \La$ we compute $\HH^* ( \La , {_1
\Lambda_{\psi}} ) = \Ext_{\Lae}^*( \La , {_1 \Lambda_{\psi}} )$ and
$\HH_* ( \La , {_1 \Lambda_{\psi}} ) = \Tor^{\Lae}_*( \La , {_1
\Lambda_{\psi}} )$, that is, the Hochschild cohomology and homology
of $\La$ with coefficients in the twisted bimodule
${_1\Lambda}_{\psi}$ (the action of $\Lae$ on ${_1\Lambda}_{\psi}$
is defined as $\lambda (\lambda_1 \otimes \lambda_2)= \lambda_1
\lambda \psi ( \lambda_2 )$). As a result we obtain several
interesting examples, both in cohomology and homology, of the
homological behavior of $\La$ as a bimodule.

\section{The Hochschild Homology}\label{sec2}

Denote by $D$ the usual $k$-dual $\Hom_k(-,k)$, and consider the map
$\phi \colon _{\La}\La \to D( \La_{\La} )$ of left $\La$-modules
defined by
$$\phi(1)( \alpha + \beta x + \gamma y + \delta yx )
\stackrel{\text{def}}{=} \delta.$$ It is easy to show that this is
an injective map and hence also an isomorphism since $\dim_k \La =
\dim_k D( \La )$, and therefore $\La$ is a Frobenius algebra by
definition. Now take any element $\lambda \in \La$, and consider the
element $\phi(1) \cdot \lambda \in D( \La)$ (we consider $D( \La )$
as a $\La$-$\La$-bimodule). As $\phi$ is surjective, there is an
element $\lambda' \in \La$ such that $\lambda' \cdot \phi(1) = \phi(
\lambda' ) = \phi(1) \cdot \lambda$, and the map $\lambda \mapsto
\lambda'$ defines a $k$-algebra automorphism $\nu^{-1} \colon \La
\to \La$ whose inverse $\nu$ is called the \emph{Nakayama
automorphism} of $\La$ (with respect to the map $\phi$).
Straightforward calculations show that $x \cdot \phi(1) = \phi(1)
\cdot (-q^{-1}x)$ and $y \cdot \phi(1) = \phi(1) \cdot (-qy)$, hence
since $x$ and $y$ generate $\La$ over $k$ we see that $\nu$ is the
degree preserving map defined by
$$x \mapsto -q^{-1}x, \hspace{3mm} y \mapsto -qy.$$

\sloppy The map $\phi$ induces a bimodule isomorphism
$_1\La_{\nu^{-1}} \simeq D( \La )$, which in turn gives an
isomorphism $_{(\nu \psi)^{-1}}\La_{\nu^{-1}} \simeq _{(\nu
\psi)^{-1}}D( \La )_1$ for any automorphism $\psi$ of $\La$.
Furthermore, since $_{(\nu \psi)^{-1}}\La_{\nu^{-1}}$ is isomorphic
to $_1\La_{\psi}$ and $_{(\nu \psi)^{-1}}D( \La )_1 = D( _1\La_{(\nu
\psi)^{-1}} )$, we get an isomorphism $_1\La_{\psi} \simeq D(
_1\La_{(\nu \psi)^{-1}} )$ of bimodules. Now from \cite[Proposition
VI.5.1]{Cartan} we get
\begin{eqnarray*}
\HH^n( \La , {_1\Lambda_{\psi}} ) & = & \Ext_{\Lae}^n( \La , {_1
\Lambda_{\psi}} ) \\
& \simeq & \Ext_{\Lae}^n \left ( \La ,D( _1\La_{(\nu \psi)^{-1}} )
\right ) \\
& \simeq & D \left ( \Tor^{\Lae}_n( \La, {_1\La_{(\nu \psi)^{-1}}})
\right ) \\
& = & D \left ( \HH_n( \La, {_1\La_{(\nu \psi)^{-1}}}) \right ),
\end{eqnarray*}
thus when computing the (dimension of the) Hochschild cohomology
group $\HH^n( \La , {_1\Lambda_{\psi}} )$ we are also computing the
Hochschild homology group $\HH_n( \La , _1\La_{(\nu \psi)^{-1}})$.
Moreover, as $\psi$ ranges over all degree preserving $k$-algebra
automorphisms of $\La$, so does $(\nu \psi)^{-1}$.

\section{The Cohomology Complex}\label{sec3}

\sloppy We start by recalling the construction of the minimal
bimodule projective resolution of $\La$ from \cite{Buchweitz}.
Define the elements
$$f^0_0 = 1, \hspace{3mm} f^1_0 = x, \hspace{3mm} f^1_1 = y, $$
$$f^n_{-1} =0= f^n_{n+1} \hspace{1mm} \text{ for each } n \geq 0,$$
and for each $n \geq 2$ define elements $\{f^n_i\}_{i=0}^n \subseteq
\underbrace{ \La \otimes_k \cdots \otimes_k \La}_{n \text{ copies}}$
inductively by
$$f^n_i = f^{n-1}_{i-1} \otimes y + q^i f^{n-1}_i \otimes x.$$
Denote by $P^n$ the $\Lae$-projective module $\bigoplus_{i=0}^n \La
\otimes_k f^n_i \otimes_k \La$, and by $\tilde{f}^n_i$ the element
$1 \otimes f^n_i \otimes 1 \in P^n$ (and $\tilde{f}^0_0 = 1 \otimes
1$). The set $\{ \tilde{f}^n_i \}_{i=0}^n$ generates $P^n$ as a
$\Lae$-module. Now define a map $\delta_n \colon P^n \to P^{n-1}$ by
$$\tilde{f}^n_i \mapsto \left [ x \tilde{f}^{n-1}_i + (-1)^n q^i
\tilde{f}^{n-1}_i x \right ] + \left [ q^{n-i} y
\tilde{f}^{n-1}_{i-1} + (-1)^n \tilde{f}^{n-1}_{i-1} y \right ].$$
It is shown in \cite{Buchweitz} that
$$( \mathbb{P}, \delta ) \colon \cdots \to P^{n+1}
\xrightarrow{\delta_{n+1}} P^n \xrightarrow{\delta_n} P^{n-1} \to
\cdots$$ is a minimal $\Lae$-projective resolution of $\La$.
Denote the direct sum of $n$ copies of ${_1 \Lambda_{\psi}}$ by
${_1 \Lambda_{\psi}^n}$, and consider its standard $k$-basis $\{
e^{n-1}_i, xe^{n-1}_i, ye^{n-1}_i, yxe^{n-1}_i \}_{i=0}^{n-1}$.
Define a map $d_n \colon {_1 \Lambda_{\psi}^n} \to {_1
\Lambda_{\psi}^{n+1}}$ by
$$\lambda e^{n-1}_i \mapsto \left [ x \lambda + (-1)^n q^i \lambda
\psi (x) \right ] e^n_i + \left [ q^{n-i-1} y \lambda + (-1)^n
\lambda \psi (y) \right ] e^n_{i+1}.$$ Applying $\Hom_{\Lae} ( -,
{_1 \Lambda_{\psi}} )$ to the resolution $( \mathbb{P}, \delta )$,
keeping in mind that $\Hom_{\Lae} ( P^n, {_1 \Lambda_{\psi}} )$
and ${_1 \Lambda_{\psi}^{n+1}}$ are isomorphic as $k$-vector
spaces, we get the commutative diagram
$$\xymatrix{
\cdots \ar[r] & \Hom_{\Lae} ( P^{n-1}, {_1 \Lambda_{\psi}} )
\ar[d]^{\wr} \ar[r]^{\delta_n^*} & \Hom_{\Lae} ( P^n, {_1
\Lambda_{\psi}} ) \ar[d]^{\wr} \ar[r]^-{\delta_{n+1}^*} & \cdots
\\ \cdots \ar[r] & {_1 \Lambda_{\psi}^n} \ar[r]^{d_n} & {_1
\Lambda_{\psi}^{n+1}} \ar[r]^{d_{n+1}} & \cdots }$$ of $k$-vector
spaces.

In order to compute $\HH^n ( \La , {_1 \Lambda_{\psi}} ) =
\Ext^n_{\Lae} ( \La, {_1 \Lambda_{\psi}} )$ for $n>0$ we compute the
cohomology $\Ker d_{n+1} / \Im d_n$ of the bottom complex in the
above commutative diagram. We do this by finding $\dim_k \Im d_n$;
once we know $\dim_k \Im d_n$, we obtain $\dim_k \Ker d_n$ (and
therefore also $\dim_k \Ker d_{n+1}$) from the equation
$$\dim_k \Ker d_n + \dim_k \Im d_n = \dim_k {_1 \Lambda_{\psi}^n}
= 4n.$$ We then have
$$\dim_k \HH^n ( \La , {_1 \Lambda_{\psi}} ) =
\dim_k \Ker d_{n+1} - \dim_k \Im d_n.$$

Now let $\psi$ be a degree preserving $k$-algebra automorphism of
$\La$. Then there are elements $\alpha_1, \alpha_2, \beta_1, \beta_2
\in k$ such that $\psi (x) = \alpha_1 x + \alpha_2 y$ and $\psi (y)
= \beta_1 x + \beta_2 y$. Since $\psi (x^2) = \psi (y^2) = \psi
(xy+qyx) =0$, we have the relations $\alpha_1 \alpha_2 = \beta_1
\beta_2 = \alpha_2 \beta_1 = 0$. If $\alpha_2 \neq 0$, then
$\alpha_1 = \beta_1 =0$, implying $x \notin \Im \psi$. Similarly, if
$\beta_1 \neq 0$, then $\alpha_2 = \beta_2 =0$, implying $y \notin
\Im \psi$. Therefore $\alpha_2 = \beta_1 = 0$, and this forces
$\alpha_1$ and $\beta_2$ to be nonzero. Thus the degree preserving
$k$-algebra automorphisms of $\La$ are precisely those defined by
$$x \mapsto \alpha x, \hspace{3mm} y \mapsto \beta y$$
for two arbitrary nonzero elements $\alpha, \beta \in k$. For such an
automorphism, the result of applying $d_n$ to the basis vectors $\{
e^{n-1}_i, xe^{n-1}_i, ye^{n-1}_i, yxe^{n-1}_i \}_{i=0}^{n-1}$ of
${_1 \Lambda_{\psi}^n}$ is
\begin{eqnarray*}
e^{n-1}_i & \mapsto & \left [ 1 + (-1)^n q^i \alpha \right ]
xe^n_i + \left [ q^{n-i-1} + (-1)^n \beta \right ] ye^n_{i+1} \\
xe^{n-1}_i & \mapsto & \left [ q^{n-i-1} + (-1)^{n+1} q \beta
\right ] yxe^n_{i+1} \\
ye^{n-1}_i & \mapsto & \left [ -q + (-1)^n q^i \alpha \right ]
yxe^n_i \\
yxe^{n-1}_i & \mapsto & 0
\end{eqnarray*}
for $0 \leq i \leq n-1$. Note that the inequality $\dim_k \Im d_n
\leq 2n+1$ always holds.

\section{The Hochschild Cohomology}\label{sec4}

We start by computing $\HH^0( \La , {_1 \Lambda_{\psi}} )$. Rather
than computing this vector space directly using the identifications
\begin{eqnarray*}
\HH^0( \La , {_1 \Lambda_{\psi}} ) &=& \{ z \in {_1 \Lambda_{\psi}}
\mid \lambda \cdot z = z \cdot \lambda \text{ for all } \lambda \in
\Lambda \} \\
&=& \{ z \in {_1 \Lambda_{\psi}} \mid \lambda z = z \psi ( \lambda )
\text{ for all } \lambda \in
\Lambda \} \\
&=& \{ z \in {_1 \Lambda_{\psi}} \mid x z = z \psi(x), y z = z
\psi(y), y x z = z \psi(y x) \} \\
&=& \{ z \in {_1 \Lambda_{\psi}} \mid x z = \alpha z x, y z = \beta
z y, y x z = \alpha \beta z y x \},
\end{eqnarray*}
we use our cohomology complex and the isomorphism $\HH^0( \La , {_1
\Lambda_{\psi}} ) \simeq \Ker d_1$. From the above we see that the
map $d_1$ is defined by
\begin{eqnarray*}
e^0_0 & \mapsto & [1- \alpha ]xe^1_0 + [1- \beta ]ye^1_1 \\
xe^0_0 & \mapsto & [1+q \beta ]yxe^1_1 \\
ye^0_0 & \mapsto & -[q+ \alpha]yxe^1_0 \\
yxe^0_0 & \mapsto & 0,
\end{eqnarray*}
and so calculation gives
$$\dim_k \HH^0( \La , {_1 \Lambda_{\psi}} ) = \left \{
\begin{array}{ll}
 3 & \text{when } \alpha = -q, \beta =-q^{-1}\\
 2 & \text{when } \alpha =1, \beta =1\\
 2 & \text{when } \alpha = -q, \beta \neq -q^{-1}\\
 2 & \text{when } \alpha \neq -q, \beta = -q^{-1}\\
 1 & \text{otherwise}
 \end{array} \right. $$
when the characteristic of $k$ is not $2$. In the characteristic $2$
case we replace $-q, -q^{-1}$ and $1$ in the above formula by $\pm
q, \pm q^{-1}$ and $\pm 1$, respectively.

Now we turn to the cohomology groups $\HH^n( \La , {_1
\Lambda_{\psi}} )$ for $n>0$. To compute the dimension of $\Im d_n$,
we distinguish between four possible cases depending on whether or
not $\alpha$ and $\beta$ belong to the set
$$\Sigma = \{ \pm q^i \}_{i \in \mathbb{Z}}.$$

\subsection{The case $\alpha, \beta \not\in \Sigma$:}

$ $

This is the easiest case; $d_n(e^{n-1}_i), d_n(xe^{n-1}_i)$ and
$d_n(ye^{n-1}_i)$ are all nonzero for $0 \leq i \leq n-1$, hence
$\dim_k \Im d_n = 2n+1$ for all $n$. Then $\dim_k \Ker d_n =2n-1$,
implying $\dim_k \Ker d_{n+1}= 2n+1$ and therefore that $\HH^n ( \La
, {_1 \Lambda_{\psi}} ) =0$ for $n > 0$.

\begin{remark}
We can relate the vanishing of cohomology to the conjecture of
Tachikawa stating that over a selfinjective ring the only finitely
generated modules having no self extensions are the projective ones.
Namely, let $M$ be a finitely generated $\La$-module such that $\La$
has no bimodule extensions by $\Hom_k(M,M)$. Then from
\cite[Corollary IX.4.4]{Cartan} we get
$$\Ext_{\La}^n(M,M) \simeq \HH^n (\La, \Hom_k (M,M) ) = 0$$
for $n>0$. Since Tachikawa's conjecture holds for $\La$ (see
\cite[Proposition 4.2]{Schulz}), the module $M$ must be projective
and therefore ($\La$ is local) isomorphic to $\La^t$ for some $t \in
\mathbb{N}$. This gives
$$\Hom_k (M,M) \simeq \left ( \La \otimes_k D(\La) \right )^{t^2}
\simeq (\Lae)^{t^2}.$$ In particular, if $\HH^n ( \La , {_1
\Lambda_{\psi}} ) =0$ for $n > 0$ then there cannot exist a
$\La$-module $M$ such that $\Hom_k (M,M)$ is isomorphic to ${_1
\Lambda_{\psi}}$, since this would imply the contradiction ${_1
\Lambda_{\psi}} \simeq (\Lae)^{t^2}$.
\end{remark}

\subsection{The case $\alpha \in \Sigma, \beta \not\in \Sigma$:}

$ $

Since $\left [ q^{n-i-1} + (-1)^n \beta \right ]$ and $\left [
q^{n-i-1} + (-1)^{n+1} q \beta \right ]$ are nonzero, we have
$d_n(e^{n-1}_i) \neq 0$ and $d_n(xe^{n-1}_i) \neq 0$ for $0 \leq i
\leq n-1$. Therefore $\dim_k \Im d_n \geq 2n$, and the problem is
now whether or not the basis vector $yxe^n_0$ belongs to $\Im d_n$.
This is the case if and only if $d_n ( ye^{n-1}_0 ) \neq 0$, that
is, if and only if
\begin{equation*}\label{condition1}
-q + (-1)^n \alpha \neq 0 \tag{C1}
\end{equation*}
holds. We now break down this case into three cases.

\subsubsection*{(i) The case $\alpha \neq \pm q$:}

$ $

Since (\ref{condition1}) holds we have $\dim_k \Im d_n = 2n+1$, and
so $\HH^n ( \La , {_1 \Lambda_{\psi}} ) =0$ for $n > 0$.

\subsubsection*{(ii) The case $\alpha = q$:}

$ $

When the characteristic of $k$ is not $2$, the condition
(\ref{condition1}) holds if and only if $n$ is odd. Therefore
$$\dim_k \Im d_n = \left \{ \begin{array}{ll}
                             2n & \text{for $n$ even} \\
                             2n+1 & \text{for $n$ odd},
                            \end{array} \right. $$
and so
$$\dim_k \Ker d_{n+1} = \left \{ \begin{array}{ll}
                             2n+1 & \text{for $n$ even} \\
                             2n+2 & \text{for $n$ odd}.
                            \end{array} \right. $$
This gives $\dim_k \HH^n ( \La , {_1 \Lambda_{\psi}} ) =1$ for
$n>0$.

When $k$ is of characteristic $2$ we see that (\ref{condition1})
never holds, hence $\dim_k \Im d_n =2n$. Then $\dim_k \Ker d_{n+1}
=2n+2$, giving $\dim_k \HH^n ( \La , {_1 \Lambda_{\psi}} ) =2$ for
$n>0$.

\subsubsection*{(iii) The case $\alpha = -q$:}

$ $

As in the previous case, we get $\dim_k \HH^n ( \La , {_1
\Lambda_{\psi}} ) =1$ for $n>0$ when $k$ is not of characteristic
$2$, and $\dim_k \HH^n ( \La , {_1 \Lambda_{\psi}} ) =2$ for $n>0$
in the characteristic $2$ case.

\begin{remark}
From this case we obtain an example showing that symmetry in the
vanishing of $\Ext$ over $\Lae$ does not hold. Namely, define $\psi$
by
$$x \mapsto q^{-1}x, \hspace{3mm} y \mapsto \beta y$$
for some element $\beta$ not contained in $\Sigma$. For $n>0$, case
(i) above gives $\Ext_{\Lae}^n ( \La, {_1 \Lambda_{\psi}} ) =0$,
whereas from case (ii) we see that $\dim_k \Ext_{\Lae}^n ( \La, {_1
\Lambda_{\psi^{-1}}} )$ is either $1$ or $2$, depending on the
characteristic of $k$. Now it is easy to see that $\Ext_{\Lae}^n (
\La, {_1 \Lambda_{\psi^{-1}}} )$ is isomorphic to $\Ext_{\Lae}^n (
{_1 \Lambda_{\psi}}, \La )$.
\end{remark}

\subsection{The case $\alpha \not\in \Sigma, \beta \in \Sigma$:}

$ $

The algebra $\La$ is isomorphic to the algebra $k \langle u,v
\rangle / ( u^2, uv+q^{-1}vu, v^2 )$ via the map
$$x \mapsto v, \hspace{3mm} y \mapsto u,$$
hence this case is symmetric to the case $\alpha \in \Sigma, \beta
\not\in \Sigma$ treated above. Namely, when $\beta \neq \pm q^{-1}$
the result is as in (i), whereas when $\beta = \pm q^{-1}$ the
result is as in (ii) and (iii).

\subsection{The case $\alpha \in \Sigma, \beta \in \Sigma$:}

$ $

The basis vectors $yxe^n_0$ and $yxe^n_n$ can only be the image of
$ye^{n-1}_0$ and $xe^{n-1}_{n-1}$, respectively, whereas for $1 \leq
i \leq n-1$ the basis vector $yxe^n_i$ can be the image of both
$ye^{n-1}_i$ and $xe^{n-1}_{i-1}$. Therefore we break this case down
into four cases, each depending on whether or not $\alpha = \pm q$
and $\beta = \pm q^{-1}$.

\subsubsection*{(i) The case $\alpha = \pm q, \beta = \pm
q^{-1}$:}

$ $

We have
$$d_n ( e^{n-1}_i ) = \left [ 1 \pm (-1)^n q^{i+1} \right ]
xe^n_i + \left [ q^{n-i-1} \pm (-1)^n q^{-1} \right ] ye^n_{i+1},$$
and since $i+1 \geq 1$ the term $\left [ 1 \pm (-1)^n q^{i+1} \right
]$ must be nonzero. Therefore $d_n ( e^{n-1}_i ) \neq 0$. Applying
$d_n$ to $xe^{n-1}_i$ and $ye^{n-1}_i$ gives $\left [ q^{n-i-1} \pm
(-1)^{n+1} \right ] yxe^n_{i+1}$ and $\left [ -q \pm (-1)^n q^{i+1}
\right ] yxe^n_i$, respectively, hence when the characteristic of
$k$ is not $2$ we get
\begin{eqnarray*}
d_n ( xe^{n-1}_i ) &=& \left \{ \begin{array}{ll}
                                 0 & \text{for $i=n-1, \beta =
                                 q^{-1}, n$ even} \\
                                 0 & \text{for $i=n-1, \beta =
                                 -q^{-1}, n$ odd} \\
                                 \neq 0 & \text{otherwise},
                                \end{array} \right. \\
d_n ( ye^{n-1}_i ) &=& \left \{ \begin{array}{ll}
                                 0 & \text{for $i=0, \alpha =q, n$
                                 even} \\
                                 0 & \text{for $i=0, \alpha =-q,
                                 n$ odd} \\
                                 \neq 0 & \text{otherwise}.
                                \end{array} \right.
\end{eqnarray*}
There are four possible pairs $( \alpha, \beta )$ to consider. If
$\alpha =q$ and $\beta =q^{-1}$, then the above gives
$$\dim_k \Im d_n = \left \{ \begin{array}{ll}
                             2n-1 & \text{for $n$ even} \\
                             2n+1 & \text{for $n$ odd},
                            \end{array} \right. $$
and therefore
$$\dim_k \Ker d_{n+1} = \left \{ \begin{array}{ll}
                             2n+1 & \text{for $n$ even} \\
                             2n+3 & \text{for $n$ odd}.
                            \end{array} \right. $$
This implies $\dim_k \HH^n ( \La , {_1 \Lambda_{\psi}} ) =2$ for
$n>0$. Similar computation gives $\dim_k \HH^n ( \La , {_1
\Lambda_{\psi}} ) =2$ for $n>0$ also for the other three possible
pairs.

When $k$ is of characteristic $2$ then
\begin{eqnarray*}
d_n ( xe^{n-1}_i ) &=& \left \{ \begin{array}{ll}
                                 0 & \text{for $i=n-1$} \\
                                 \neq 0 & \text{otherwise},
                                \end{array} \right. \\
d_n ( ye^{n-1}_i ) &=& \left \{ \begin{array}{ll}
                                 0 & \text{for $i=0$} \\
                                 \neq 0 & \text{otherwise},
                                \end{array} \right.
\end{eqnarray*}
and so $\dim_k \Im d_n = 2n-1$ for all $n>0$. Consequently $\dim_k
\HH^n ( \La , {_1 \Lambda_{\psi}} ) =4$ for $n>0$.

\begin{remark}
Note that, in this particular case, we have computed the (dimension
of the) Hochschild homology $\HH_* (\La) = \Tor^{\Lae}_* ( \La, \La
)$ of $\La$. Namely, it follows from Section \ref{sec2} that for
each $n>0$ the $k$-vector spaces $\HH^n( \La, _1\La_{\nu^{-1}} )$
and $D\left ( \HH_n( \La ) \right )$ are isomorphic, where $\nu$ is
the Nakayama automorphism
$$x \mapsto -q^{-1}x, \hspace{3mm} y \mapsto -qy$$
of $\La$. Then $\nu^{-1}$ is defined by
$$x \mapsto -qx, \hspace{3mm} y \mapsto -q^{-1}y,$$
hence ${_1 \Lambda_{\nu^{-1}}}$ is precisely the sort of bimodule we
have just considered in terms of Hochschild cohomology. Consequently
$\dim_k \HH_n (\La) =2$ for $n>0$ when the characteristic of $k$ is
not $2$, whereas $\dim_k \HH_n (\La) =4$ for $n>0$ in the
characteristic $2$ case.
\end{remark}

\subsubsection*{(ii) The case $\alpha = \pm q, \beta \neq \pm
q^{-1}$:}

$ $

As in (i) the element $d_n ( e^{n-1}_i )$ is nonzero for $0 \leq i
\leq n-1$. Moreover, since $\beta \neq \pm q^{-1}$ the basis element
$yxe^n_n$ always lies in $\Im d_n$, as does the basis elements
$yxe^n_i$ for $1 \leq i \leq n-1$. Therefore $\dim \Im d_n \geq 2n$,
and the question is whether or not $yxe^n_0$ belongs to $\Im d_n$,
i.e.\ whether or not $d_n( ye^{n-1}_0)$ is nonzero.

When the characteristic of $k$ is not $2$, then from (i) we see that
$$d_n ( ye^{n-1}_0 ) = \left \{ \begin{array}{ll}
                                 0 & \text{for $\alpha =q, n$
                                 even} \\
                                 0 & \text{for $\alpha =-q,
                                 n$ odd} \\
                                 \neq 0 & \text{otherwise},
                                \end{array} \right. $$
and computation gives $\dim_k \HH^n ( \La , {_1 \Lambda_{\psi}} )
=1$ for $n>0$. However, when the characteristic of $k$ is $2$ then
$d_n( ye^{n-1}_0)=0$, giving $\dim_k \Im d_n = 2n$ and consequently
$\dim_k \HH^n ( \La , {_1 \Lambda_{\psi}} ) = 2$ for $n>0$.

\subsubsection*{(iii) The case $\alpha \neq \pm q, \beta = \pm
q^{-1}$:}

$ $

Using the isomorphism $\La \simeq k \langle u,v \rangle / ( u^2,
uv+q^{-1}vu, v^2 )$ from the case $\alpha \not\in \Sigma, \beta \in
\Sigma$, we see that the present case is symmetric to the case (ii)
above. Thus when the characteristic of $k$ is not $2$ then $\dim_k
\HH^n ( \La , {_1 \Lambda_{\psi}} ) =1$ for $n>0$, whereas $\dim_k
\HH^n ( \La , {_1 \Lambda_{\psi}} ) = 2$ for $n>0$ in the
characteristic $2$ case.

\subsubsection*{(iv) The case $\alpha \neq \pm q, \beta \neq \pm
q^{-1}$:}

$ $

We now have $\alpha = \pm q^s$ and $\beta = \pm q^t$ where $s \in
\mathbb{Z} \setminus \{ 1 \}$ and $t \in \mathbb{Z} \setminus \{
-1 \}$, and therefore the basis elements $yxe^n_0$ and $yxe^n_n$
both lie in the image of $d_n$. To compute the dimension of $\Im
d_n$, we must find out when $d_n ( e^{n-1}_i ) =0$ for some $0
\leq i \leq n-1$ and when $yxe^n_i \not\in \Im d_n$ for some $1
\leq i \leq n-1$. We have
\begin{eqnarray*}
& \left \{ d_n ( e^{n-1}_i ) =0 \text{ for some } 0 \leq i \leq
n-1 \right \} & \\
(*) & \Updownarrow & \\
& \left \{ 1+ (-1)^n q^i \alpha =0 \text{ and } q^{n-i-1} + (-1)^n
\beta =0 \right \} &
\end{eqnarray*}
and
\begin{eqnarray*}
& \left \{ yxe^n_i \not\in \Im d_n \text{ for some } 1 \leq i \leq
n-1 \right \} & \\
(**) & \Updownarrow & \\
& \left \{ -q +(-1)^nq^i \alpha =0 \text{ and } q^{n-i} +(-1)^{n+1}q
\beta =0 \right \}, &
\end{eqnarray*}
and when the characteristic of $k$ is not $2$ this happens
precisely when we have the following:
\begin{equation*}\label{condition2}
s \leq 0, \hspace{3mm} t \geq 0, \hspace{3mm} \alpha =
(-1)^{t-s}q^s, \hspace{3mm} \beta = (-1)^{t-s}q^t \tag{C2}.
\end{equation*}
In the characteristic $2$ case we may relax this condition; in
this case ($*$) and ($**$) occur precisely when we have the
following:
\begin{equation*}\label{condition3}
s \leq 0, \hspace{3mm} t \geq 0 \tag{C3}.
\end{equation*}
However, both when the characteristic of $k$ is not $2$ and
(\ref{condition2}) holds, and in the characteristic $2$ case when
(\ref{condition3}) holds, we see that ($*$) occurs for $n =
t-s+1$, whereas ($**$) occurs for $n=t-s+2$.

Therefore, when the characteristic of $k$ is not $2$ the dimension
of $\Im d_n$ is given by
$$\dim_k \Im d_n = \left \{ \begin{array}{ll}
                             2n & \text{when (\ref{condition2})
                             holds and } n=t-s+1 \\
                             2n & \text{when (\ref{condition2})
                             holds and } n=t-s+2 \\
                             2n+1 & \text{otherwise},
                             \end{array} \right. $$
implying the dimension of $\Ker d_{n+1}$ is given by
$$\dim_k \Ker d_{n+1} = \left \{ \begin{array}{ll}
                                  2n+2 & \text{when
                                  (\ref{condition2}) holds and }
                                  n=t-s \\
                                  2n+2 & \text{when
                                  (\ref{condition2}) holds and }
                                  n=t-s+1 \\
                                  2n+1 & \text{otherwise}.
                                  \end{array} \right. $$
Consequently , the dimension of $\HH^n ( \La , {_1 \Lambda_{\psi}}
)$ for $n>0$ is given by
$$\dim_k \HH^n ( \La , {_1 \Lambda_{\psi}} ) = \left \{
\begin{array}{ll}
1 & \text{when (\ref{condition2}) holds and } n=t-s \\
2 & \text{when (\ref{condition2}) holds and } n=t-s+1 \\
1 & \text{when (\ref{condition2}) holds and } n=t-s+2 \\
0 & \text{otherwise}.
\end{array} \right. $$
When the characteristic of $k$ is $2$, we obtain the exact same
formulas, but with (\ref{condition2}) replaced by
(\ref{condition3}).

\begin{remark}
(i) In the case considered above, we see that the cohomology is zero
except possibly in three degrees, depending on what conditions
$s,t,\alpha$ and $\beta$ satisfy. As a consequence, we construct a
counterexample to the following conjecture by Auslander (see
\cite[page 815]{Auslander}): if $M$ is a finitely generated module
over an Artin algebra $\Gamma$, then there exists a number $n_M$
such that for any finitely generated module $N$ we have
$$\Ext_{\Gamma}^i (M,N)=0 \text{ for } i \gg 0 \Rightarrow
\Ext_{\Gamma}^i (M,N)=0 \text{ for } i \geq n_M.$$ The first
counterexample to this conjecture appeared in \cite{Jorgensen},
where the algebra considered was a finite dimensional commutative
Noetherian local Gorenstein algebra. A counterexample over our
algebra $\La = k \langle x,y \rangle / ( x^2, xy+qyx, y^2 )$ was
given in \cite{Smalo}.

As for a counterexample using Hochschild cohomology, define for
each natural number $t$ an automorphism $\psi \colon \La \to \La$
by $x \mapsto q^{-t}x$ and $y \mapsto q^ty$, and denote the
bimodule ${_1 \Lambda_{\psi}}$ by $M_t$. Then condition
(\ref{condition2})/(\ref{condition3}) is satisfied (with
$s=-t,\alpha=q^{-t}$ and $\beta=q^t$), and so
$$\Ext_{\Lae}^n ( \La, M_t ) = \left \{ \begin{array}{ll}
                                        \neq 0 & \text{for }
                                        n=2t+2\\
                                        0 & \text{for } n>2t+2.
                                        \end{array} \right. $$

(ii) Even though the above conjecture of Auslander fails to hold in
general, Auslander himself proved (unpublished, see \cite[page
815]{Auslander}) that if the conjecture holds for the enveloping
algebra $\Gamma^{\e}$ of a finite dimensional algebra $\Gamma$ over
a field, then the finitistic dimension
$$\sup \{ \Pd_{\Gamma} X \mid X \text{ finitely generated
$\Gamma$-module with } \Pd_{\Gamma} X < \infty \}$$ of $\Gamma$ is
finite. In view of the above remark, we see that the converse to
this result does not hold; our algebra $\La$, being selfinjective,
trivially has finite finitistic dimension, whereas the conjecture
of Auslander does not hold for $\Lae$.

(iii) The computation of the Hochschild cohomology $\HH^* ( \La ) =
\Ext_{\Lae}^* ( \La, \La )$ of $\La$ is covered by the last of the
above cases. When $\psi$ is the identity automorphism we have
$s=t=0$, and the condition (\ref{condition2})/(\ref{condition3}) is
satisfied. This gives
$$\dim_k \HH^n ( \La ) = \left \{ \begin{array}{ll}
                                    2 & \text{for } n=1 \\
                                    1 & \text{for } n=2 \\
                                    0 & \text{for } n \geq 3,
                                  \end{array} \right. $$
and so our algebra $\La$ is a counterexample to the following
question raised by Happel in \cite{Happel}: if the Hochschild
cohomology groups of a finite dimensional algebra vanish in high
degrees, does the algebra have finite global dimension? The
counterexample above first appeared in \cite{Buchweitz}, where it
was shown that the generating function $\sum_{n=0}^{\infty} \HH^n (
\La )t^n$ of $\HH^n ( \La )$ is $2+2t+t^2$.

The converse to the question of Happel is always true when the
algebra modulo its Jacobson radical is separable over the ground
field. More specifically, if $\Gamma$ is a finite dimensional
algebra over a field $K$ with Jacobson radical $\ra$, and the
semisimple algebra $\Gamma / \ra$ is separable over $K$, then by
\cite[\S 3]{Eilenberg} the implication
$$\gldim \Gamma < \infty \Rightarrow \Pd_{\Gamma^{\e}} \Gamma <
\infty$$ holds. In particular the Hochschild cohomology groups
$\HH^n ( \Gamma ) = \Ext_{\Gamma^{\e}}^n ( \Gamma, \Gamma )$ and the
homology groups $\HH_n ( \Gamma ) = \Tor^{\Gamma^{\e}}_n ( \Gamma,
\Gamma )$ of $\Gamma$ vanish for $n \gg 0$ when the global dimension
of $\Gamma$ is finite. It is not known whether the vanishing of the
Hochschild homology groups in high degrees for a finite dimensional
algebra implies the global dimension of the algebra is finite.

(iv) Because of the equality $\dim_k \HH^n( \La , {_1\Lambda_{\psi}}
) = \dim_k \HH_n( \La, {_1\La_{(\nu \psi)^{-1}}})$ which follows
from Section \ref{sec2}, the somewhat strange behavior in Hochschild
cohomology revealed in the last case considered above can also be
transferred to Hochschild homology. When the automorphism $\psi$ is
given by $\psi (x) = \pm q^s x$ and $\psi (y) = \pm q^t y$, then the
automorphism $\theta \stackrel{\text{def}}{=} (\nu \psi )^{-1}$,
where $\nu$ is the Nakayama automorphism, is given by
$$x \mapsto \mp q^{1-s}x, \hspace{3mm} y \mapsto \mp q^{-(t+1)}y.$$
Thus for such an automorphism $\theta$, when $s \in \mathbb{Z}
\setminus \{ 1 \}$ and $t \in \mathbb{Z} \setminus \{ -1 \}$ we get
$$\dim_k \HH_n ( \La , {_1 \Lambda_{\theta}} ) = \left \{
\begin{array}{ll}
1 & \text{when (\ref{condition2})/(\ref{condition3}) holds and }
n=t-s \\
2 & \text{when (\ref{condition2})/(\ref{condition3}) holds and }
n=t-s+1 \\
1 & \text{when (\ref{condition2})/(\ref{condition3}) holds and }
n=t-s+2 \\ 0 & \text{otherwise}.
\end{array} \right. $$
when $n>0$. In these formulas condition (\ref{condition2}) applies
when the characteristic of $k$ is not $2$, and condition
(\ref{condition3}) applies in the characteristic $2$ case.
\end{remark}

\section*{Acknowledgement}

I would like to thank my supervisor {\O}yvind Solberg for valuable
suggestions and comments.

\end{document}